\documentclass{article}
\usepackage{amsthm,amsfonts, amsbsy, amssymb,amsmath,graphicx}
\usepackage[russian]{babel}

\newtheorem{theorem}{Òåîðåìà}

\begin{document}

%\begin{tabular}{lr}
%ÓÄÊ 517.956.3&~~~~~~~~~~~~~~~~~~~~~~ÌÀÒÅÌÀÒÈÊÀ
%\end{tabular}
%\vspace{5mm}
\title{SAMUELSON'S WEBS}
\author{@ 2009 ã., V.~V.~Goldberg and V.~V.~Lychagin\footnote{\textit{New Jersey Institute of Technology, Newark, New Jersey, United States of America; Institute of Control Problems, Moscow, Russia; email:
vladislav.goldberg@gmail.com, lychagin@yahoo.com}}}
\date{}
\maketitle

\begin{abstract}
In the present paper we define Samuelson's webs
and their rank. The main result of the paper is the proof that the rank of the Samuelson webs does not exceed 6, as well as finding the conditions under which this rank is maximal for the general Samuelson webs as well as for their singular cases.
\end{abstract}

\section{Introduction}
Application of web theory in economics has its origins in the works of Gerard Debreu and Paul A. Samuelson (Nobel Prize winners in economics in 1983 and 1972).

Debreu showed that the existence of a function providing a preference ordering is equivalent to the triviality of a certain $3$-web. In Samuelson's theory the problem of profit maximization can be formulated in terms of web theory (see, for example, \cite{CR09}).

In this paper we propose interpreting the area condition of Maxwell--Samuelson in terms of webs as a  quadratic relation on the differential forms defining a planar $4$-web. This makes transparent the relationship between Maxwell--Samuelson's condition and Abel's equations. We use this observation to derive a system of differential equations, which we call Samuelson's equations, and which to a large extent are similar to the Abel equations. As for Abel's equations, we introduce the notion of the rank of Samuelson's web ($S$-rank) which coincides with the dimension of the solution space of the system of Samuelson's equations.

The main result of this paper is the proof
that the rank of the Samuelson webs does not exceed 6, as well as finding the conditions under which this rank is maximal. Our approach is constructive and, in particular, it contains a procedure for finding the rank of an arbitrary $S$-web (cf. \cite{GL09}).

\section{Samuelson's webs}
Let $M = <\omega_1, \omega_2,\omega_3,\omega_4>$ be a  $4$-web in the plane. We say that this $4$-web is  a  \textit{Samuelson's web}, if some differential forms
$\omega_1, \omega_2,\omega_3$ and $\omega_4$ defining the $4$-web satisfy the following exterior quadratic relation:
\begin{equation}
\omega_{3}\wedge \omega_{1}+\omega_{4}\wedge \omega_{2}=0.
\label{S-condition}
\end{equation}

In what follows, for brevity we shall call Samuelson's webs $S$\emph{-webs}.

Symplectic and contact geometry provide examples of $S$-webs. Let $\mathbb{R}^4$ be a four-dimensional symplectic manifold with a structure form $dy_1 \wedge dx_1 + dy_2 \wedge dx_2$, and let
$M^2 \subset \mathbb{R}^4$ be a Lagrangian surface on which the differentials of any pair of the coordinate functions $(x_1, x_2), (y_1, y_2), (x_i, y_j)$ are linearly independent. Then the $4$-web on this surface defined by the level curves of these functions is an $S$-web.

We shall call this web the \emph{coordinate $S$-web} on the Lagrangian surface.

In a similar way an $S$-web arises on Legendrian manifolds  $M^2 \subset \mathbb{R}^5$, lying in contact manifolds $(\mathbb{R}^5, dz - ydx)$.

In the definition of $S$-webs the differential 1-forms $\omega_i$ are defined up to multiplication by functions $\lambda_i$ which do not vanish at any point:
\begin{equation}
\omega_i \rightarrow \lambda_i \omega_i,
\end{equation}
and which satisfy the condition
\begin{equation*}
\lambda_{3}\lambda_{1}=\lambda_{4}\lambda_{2}.
\label{lambda-condition}
\end{equation*}

This condition
allows us to make the following normalization of the $4$-web. First, we can choose the factors $\lambda_{1}, \lambda_{2}$ è $\lambda_{3}$
in such a way that
\begin{equation}
\omega_{3} + \omega_{1} + \omega_{2}=0.
\label{1st normalization}
\end{equation}
With this choice of the differential forms $\omega_{1}, \omega_{2}$ and $\omega_{3}$, one can show that the factors $\lambda_i$ in $\omega_i$ in relation (\ref{lambda-condition}) must be equal:
\begin{equation*}
\lambda_{1} = \lambda_{2} = \lambda_{3} = \lambda_{4}=\lambda.
\end{equation*}

With normalization (\ref{1st normalization}), condition  (\ref{S-condition}) takes the form
$
(\omega_{4} + \omega_{1}) \wedge \omega_{2} = 0.
$
It follows that
\begin{equation*}
\omega_{4} + \omega_{1} + b\omega_{2}=0,
\end{equation*}
where $b$ gives the basis invariant of the 4-web (see \cite{GL06}).

More precisely, the value of the function $b^{-1}(a)$ at the point $a \in \mathbb{R}^2$ is the cross-ratio of the points $[\omega_{1,a}, \omega_{2,a}, \omega_{3,a}, \omega_{4,a}]$ on the projective straight line $\mathbb{P} (T_{a}^*(\mathbb{R}^2) )$.

Secondly, we choose the factor $\lambda$  in such a way that $\omega_3 = df$ for some function $f$. Denote by $x$ and $y$ such functions which satisfy the conditions
$$
   \omega_1 \wedge dx = 0, \;\; \omega_2 \wedge dy = 0,
$$
but where $dx \wedge dy \neq 0$.

 These functions $x$ and $y$ can be chosen as coordinates, and equation (\ref{1st normalization}) gives the following relations:
\begin{equation*}
\omega_3 = df, \;\; \omega_1 = - f_x dx, \;\; \omega_2 = - f_y dy.
\end{equation*}

\section{Structure equations}
As in \cite{GL06}, denote by $\gamma$ such a 1-form that
$$
d\omega_i = \omega_i \wedge \gamma, \;\; i = 1, 2,  3.
$$
This form determines the Chern connection in the plane, and the first three web foliations are geodesic with respect to this connection.  The curvature form $d\gamma$ of this connection  is a differential invariant of the $3$-web $<\omega_1, \omega_2,\omega_3>$.

Denote by $\{\partial_1, \partial_2\}$ the basis of vector fields dual to the basis $\{\omega_1, \omega_2\}$:
$
<\omega_i, \partial_j> = \delta_{ij}, \;\; i. j = 1, 2.
$

Then for any function $h$ we have
\begin{equation}
dh = h_1 \omega_1 + h_2 \omega_2, \label{dh}
\end{equation}
where $h_1 = \partial_1 (h)$ and $h_2 = \partial_2 (h)$.

Differentiating relation (\ref{dh}), we find that
\begin{equation}
[\partial_1, \partial_2] = H (\partial_2 - \partial_1),
\label{commutator}
\end{equation}
where
\begin{equation*}
H = \displaystyle \frac{f_{xy}}{f_x f_y}.
\end{equation*}

In the coordinates $(x, y)$, the vector fields
$\partial_1$ and $\partial_2$ have the following form:
$$
\partial_1 = - f_x^{-1} \partial_x, \;\; \partial_2 = - f_y^{-1} \partial_y.
$$

In what follows we shall use the notation:
$
h_i = \partial_i (h), \;\; h_{ij} = \partial_i \partial_j (h), \;\; \textrm{etc}.
$

\section{Samuelson's equations}
The Maxwell--Samuelson area condition means that an  $S$-web can be realized as indicated in our main example.

Precisely this relation requires that the web $W_4$ is equivalent to a coordinate 4-web on a certain Lagrangian surface.

More precisely, it is required that there exist positive factors $s_1, s_2, t_1$ and $t_2$ such that the forms $s_1 \omega_1, s_2 \omega_2, t_1\omega_3$ and $t_2 \omega_4$ satisfy relation (\ref{S-condition}) and are closed.

These requirements imply the following relations:
\begin{equation}
\renewcommand{\arraystretch}{1.5}
\left\{
\begin{array}{ll}
d(s_1 \omega_1) = d(s_2 \omega_2) = d(t_1\omega_3)
= d(t_2 \omega_4) = 0, \\
s_1 t_1 = s_2 t_2.
 \label{S-equation}
\end{array}%
\right.
\renewcommand{\arraystretch}{1}
\end{equation}

We find now the explicit form of these equations.

We have
$$
d (s_1 \omega_1) = (H s_1 - s_{1,2}) \omega_1 \wedge \omega_2 = 0.
$$
It follows that
$$
s_{1,2} = H s_1.
$$

Similarly, we get
$$
s_{2,1} = H s_2.
$$

For the third equation of system (\ref{S-equation}) we have
$$
d (t_1 \omega_3) = (t_{1,2} - t_{1,1}) \omega_1\wedge \omega_2 = 0,
$$
whence it follows that
$$
t_{1,2} - t_{1,1} = 0.
$$
Similarly, we have
$$
\renewcommand{\arraystretch}{1.5}
\begin{array}{ll}
d (t_2 \omega_4) =
(t_{2,2} - bt_{2,1} - t_2(b_1 - (b-1)H))
\omega_1\wedge \omega_2 = 0,
\end{array}%
\renewcommand{\arraystretch}{1}
$$	
whence it follows that
$$
t_{2,2} - bt_{2,1} - t_2(b_1 - (b-1)H) = 0.
$$

We define new functions $\sigma_i$ and $\tau_i, \, i = 1, 2$ by the following formulas:
\begin{equation*}
\sigma_i = \log s_i, \;\; \tau_i = \log t_i; \; i = 1, 2.
\end{equation*}
Then equations (\ref{S-equation}) take the following form:
\begin{equation}
\renewcommand{\arraystretch}{1.5}
\left\{
\begin{array}{ll}
\sigma_{1,2} = H, \;\; \sigma_{2,1} = H, \\
\tau_{1,2} - \tau_{1,1} = 0, & \\
 b \tau_{2,1} - \tau_{2,2} = (b - 1)H - b_1.
 \label{S-equations}
\end{array}%
\right.
\renewcommand{\arraystretch}{1}
\end{equation}

In addition, the second equation of system (\ref{S-equation}) shows that
\begin{equation}
\sigma_1 + \tau_1 = \sigma_2 + \tau_2. \label{sigma+tau}
\end{equation}

Using the last equation of system (\ref{S-equations}) and (\ref{commutator}) and representing $\tau_2$ from
(\ref{sigma+tau}) in the form  $\tau_2 = \sigma_1 + \tau_1 - \sigma_2$, we obtain the final form of equations (\ref{S-equations}):
\begin{equation}
\renewcommand{\arraystretch}{1.5}
\left\{
\begin{array}{ll}
\sigma_{1,2} = H, \;\; \sigma_{2,1} = H, \\
\tau_{1,2} - \tau_{1,1} = 0, & \\
 b \sigma_{1,1} + (b-1) \tau_{1,2}+ \sigma_{2,2} = 2bH - b_1.
 \label{final S-equations}
\end{array}%
\right.
\renewcommand{\arraystretch}{1}
\end{equation}
We shall call system (\ref{final S-equations}) \emph{Samuelson's equations}, and the dimension of the solution space of this system we shall call the \emph{rank of the $S$-web}.

Let ${\cal E}_1 \subset \mathbf{J}^1 (\pi)$ be a representaion of system $(\ref{final S-equations})$ of Samuelson's equations as a submanifold in the space of 1-jets. Here we denote by $\pi: \mathbb{R}^3 \times \mathbb{R}^2 \rightarrow \mathbb{R}^2$ the trivial fiber bundle, where $\pi: (\sigma_1, \sigma_2, \tau_1, x, y) \mapsto (x, y)$
Then one can easily check that  $\mathrm{codim} \;{\cal E}_1 = 4$, and the symbol $g_1 \subset T^* \otimes \pi$ of this system is of dimension 2.

Let further ${\cal E}_2 = {\cal E}_2^{(1)} \subset
\mathbf{J}^2 (\pi)$ be the first prolongation of
the system of Samuelson's equations.

Then direct computation shows that
$\mathrm{codim} \;{\cal E}_2 = 12$, and that the symbol $g_2$ of this system, which is the first
prolongation of the symbol $g_1$, $g_2 = g_1^{(1)},$ has the dimension 1.

For the second prolongation
${\cal E}_3 = {\cal E}_1^{(2)} \subset
\mathbf{J}^3 (\pi)$ we have respectively $\mathrm{codim}\; {\cal E}_3 = 24$, and $\dim g_1^{(2)} = 0$.

In other words, the prolongations of the system of the Samuelson's equations form the following sequence of fibrations:
$$
\mathbb{R}^2 \overset{\pi}{\leftarrow}  \mathbf{J}^0 (\pi) \overset{g_1}{\leftarrow} {\cal E}_1 \overset{g_2}{\leftarrow} {\cal E}_2 \overset{0}{\leftarrow}  {\cal E}_3.
$$
Note that the condition $g_1^{(2)} = 0$ means that the projection $\pi_{3,2}: {\cal E}_3 \rightarrow {\cal E}_2$ is a diffeomorphism, and thus the second prolongation ${\cal E}_3$ defines a certain natural connection in the fiber bundle $\pi_2: {\cal E}_2 \rightarrow \mathbb{R}^2$.

\begin{theorem}
System $(\ref{final S-equations})$ of the differential Samuelson's equations is a finite type system, and the rank of the $S$-web does not exceed $\dim \pi + \dim g_1 + \dim g_1^{(1)} = 6$.
\end{theorem}

The condition that the rank of an $S$-web is maximal means that the connection indicated above is trivial, or that all obstructions to integration of  system (\ref{final S-equations}) are absent.

As in \cite{GL06}, these conditions can be calculated using the multibracket \cite{KL06}. However, below we give an alternative method which allows us to indicate not only the condition for maximum rank of the $S$-web but also to calculate this rank in concrete cases.

\section{Calculation of the rank of $S$-webs}
Consider the first equation $\sigma_{1,2} = H$ of the system of Samuelson's equations. Taking into account relation (\ref{commutator}), we can rewrite this equation in the form
\begin{equation*}
\sigma_{1,y} = - \frac{f_{xy}}{f_x f_y},
\end{equation*}
whence
$$\sigma_1 = - \log|f_x| + s_1 (x).$$
Similarly,

$$\sigma_2 = - \log |f_y| + s_2 (y).$$
The solutions of the third equation
$$
(\partial_1 - \partial_2) (\tau_1) = 0
$$
of system (\ref{final S-equations}) have the form
$$\tau_1 = w (f),$$
since the function $f$ is the first integral
of the vector field $\partial_1 - \partial_2$.

Substituting $\sigma_1, \sigma_2$, and $\tau_1$ in the last equation of system (\ref{final S-equations}), we find that
\begin{equation}
bs'_1 + s'_2 = B,
\label{s'}
\end{equation}
where
$$
B = 2bH-b_1+b\frac{f_{xx}}{f_x} + \frac{f_{yy}}{f_y} +(b-1) w'.
$$
Taking the first derivative with respect to $x$ and the second mixed derivative with respect to $x$ and $y$ of equation (\ref{s'}), we find that
\begin{equation*}
\renewcommand{\arraystretch}{1.5}
\left\{
\begin{array}{ll}
b_x s'_1 + bs''_1 = B_x, \\
b_{xy} s'_1 + b_y s''_1 = B_{xy}.
\end{array}%
\right.
\end{equation*}
Denote the determinant of this system considered as a system of linear equations with respect to the functions $s'_1$ è $s''_1$ by $\Delta$:
$$
\Delta = b_x b_y - bb_{xy}.
$$
Then if $\Delta \neq 0$, the solution of this system has the form
\begin{equation}
s'_1 = \frac{b_y B_x - b B_{xy}}{\Delta}, \label{s'1}
\end{equation}
\begin{equation}
s''_1 = \frac{{b_x} B_{xy} - b_{xy} B_{x}}{\Delta},\\
\label{s''1}
\end{equation}

Note that the function $s'_1$ depends on the variable $x$ only. Thus in order to satisfy relations (\ref{s'1}) and (\ref{s''1}), it is necessary and sufficient that the following conditions hold:
\begin{equation}
J_1 = \biggl[\frac{b_y B_x - b B_{xy}}{\Delta}\biggr]_y=0, \label{J1}
\end{equation}
\begin{equation}
J_2 = \frac{{b_x} B_{xy} - b_{xy} B_{x}}{\Delta}-\biggl[\frac{b_y B_x- b B_{xy}}{\Delta}\biggr]_x=0.
\label{J2}
\end{equation}

If these conditions are satisfied, then
equations (\ref{s'1}) and (\ref{s'}) determine the functions $s'_1 (x)$ and $s'_2 (y)$, and consequently the functions $s_1 (x)$ and $s_2 (y)$, up to additive constants.

Thus in the case when relations (\ref{J1}) and (\ref{J2}) are satisfied, the solution space of the system of equations (\ref{s'1}) and (\ref{s'}) is two-dimensional.

If at least one of equations (\ref{J1}) or
(\ref{J2}) is not satisfied and $\Delta = 0$, then system (\ref{final S-equations}) of Samuelson's equations does not have a solution.

Note also that after the substitution $b = C e^h$ the condition $\Delta = 0$ becomes the condition $h_{xy} = 0$, and thus implies that the function $b$ is a product of functions of $x$ and $y$: $b (x, y) = b_1 (x) b_2 (y)$. We shall consider this case separately.

Now consider now conditions (\ref{J1}) and
(\ref{J2}) as differential equations with respect to the function $w$. Then condition (\ref{J1}) leads to the equation
\begin{equation}
  \frac{b (1-b)}{\Delta} f_x f_y^2 w^{(4)} + T_3 w^{(3)} + T_2 w'' + T_1 w' + T_0=0.
\label{DE1 on w}
\end{equation}
Similarly from condition (\ref{J2}) we get the equation
\begin{equation}
  \frac{b (1-b)}{\Delta} f_x^2 f_y w^{(4)} + S_3 w^{(3)} + S_2 w'' + S_1 w' + S_0=0,
\label{DE2 on w}
\end{equation}
where the coefficients $T_0, T_1, T_2, T_3$ and
$S_0, S_1, S_2, S_3$ are expressed in terms of the jets of the functions $f$ and $b$ of orders five and six, respectively.

We set
$$
K_i = \frac{\Delta T_i}{b (1-b) f_x^2 f_y}, \;\; L_i = \frac{\Delta S_i}{b (1-b)f_x f_y^2}.
$$
Then equations (\ref{DE1 on w}) and (\ref{DE2 on w}) take the form
\begin{equation}
\renewcommand{\arraystretch}{1.5}
\left\{
\begin{array}{ll}
    w^{(4)} + K_3 w^{(3)} + K_2 w'' + K_1 w' + K_0=0, \\
    w^{(4)} + L_3 w^{(3)} + L_2 w'' + L_1 w' + L_0=0.
 \label{system1 for w}
\end{array}%
\right.
\end{equation}

Let us set $\delta = \partial_1 - \partial_2$. Then $\delta (w^{(i)}) = 0$, and applying the differentiation $\delta$ to (\ref{system1 for w}), we find two sequences of equations, respectively:
\begin{equation}
\renewcommand{\arraystretch}{1.5}
\left\{
\begin{array}{ll}
    \delta^{i} (K_3) w^{(3)} + \delta^{i} (K_2) w'' + \delta^{i} (K_1) w' + \delta^{i} (K_0)=0, \\
    \delta^{i} (L_3) w^{(3)} + \delta^{i} (L_2) w'' + \delta^{i} (L_1 w') + \delta^{i} (L_0)=0,
 \label{system2 for w}
\end{array}%
\right.
\end{equation}
where $i = 1, 2$.

Note that the maximal dimension of the solution space of
system (\ref{system1 for w}) equals 4.

In order to get the four-dimensional solution space,
it is necessary and sufficient that the following conditions
be satisfied:
$$
K_i = L_i, \;\; \delta (K_i) = 0, \;\; i = 0, 1, 2, 3.
$$

If the above conditions do not hold, then the dimension
of the solution space (if solutions exist) of system (\ref{system1 for w}) may vary from $-1$ to 3 depending on the behavior of the functions
$\delta (K_i)$ and $\delta (L_i)$.

\begin{theorem} Suppose that the basic invariant $b$ of an $S$-web satisfies the condition $\Delta \neq 0$. Then this $S$-web is of maximum rank $6$ if and only if
$
K_i = L_i, \; \delta (K_i) = 0
$
for all $i = 0, 1, 2, 3$.
\end{theorem}

Note that if the system of equations (\ref{system1 for w}), (\ref{J1}) and (\ref{J2}) has solutions, then the corresponding $S$-web is isomorphic to the coordinate  $S$-web on the Lagrangian surface. However, if the system indicated above does not have solutions, then the corresponding $S$-web is not isomorphic to any coordinate  $S$-web on the Lagrangian surface.

In conclusion we consider the case of singular $S$-webs, i.e., $S$-webs, satisfying the condition $\Delta = 0$.

In this case $b (x,y) = b_1 (x) b_2 (x)$, and equation (\ref{s'}) can be written in the form
\begin{equation}
b_1 s'_1 - b_2^{-1} s'_2 =b_2^{-1} B.
\label{s'bis}
\end{equation}
The condition of solvability of equation (\ref{s'bis}) is the single condition
\begin{equation}
J_3 - (b_2^{-1} B)_{xy}=0
\label{J3}
\end{equation}

Note that if this condition is satisfied, then equation (\ref{s'bis}) has a three-dimensional solution space.

In fact, differentiating equation (\ref{s'bis}) with respect to $x$, we obtain the equation $(b_1 s'_1)' = (b_2^{-1} B)_x$, which has a two-dimensional solution space.

Then, given $s_1 (x)$, equation (\ref{s'bis}) is a first-order differential equation with respect to $s_2 (y)$.

Equation (\ref{J3}) considered as equation with respect to the function $w$ has the form
\begin{equation*}
b_2^{-1} (b-1) f_x f_y w^{(3)} + r_2 w'' + r_1 w' + r_0=0,
%    \label{system3 for w}
\end{equation*}
or
$$
w^{(3)} + R_2 w'' + R_1 w' + R_0=0,
$$
where
$$
R_i = \displaystyle\frac{r_i}{(b-1) f_x f_y}, \;\; i = 0, 1, 2.
$$

Note that the coefficients $R_i$ of this equation depends on the fourth jet of the web function $f$.

\begin{theorem} Suppose that the basic invariant $b$ of an $S$-web is decomposable, i.e.,  $b (x,y) = b_1(x) b_2 (y)$. Then this $S$-web is of maximum rank $6$ if and only if $\delta (R_i) = 0, \; i = 0, 1, 2$.
\end{theorem}

\begin{center}
ACKNOWLEDGEMENTS
\end{center}
 The authors are grateful to Professors J. B. Cooper and T. Russell for fruitful discussions pertaining to the mathematical part of Samuelson's theory.

\end{document}